\documentclass[times,final]{elsarticle}

\usepackage{jcomp}
\usepackage{framed,multirow}

\usepackage{amssymb}
\usepackage{amsmath}
\usepackage{latexsym}

\usepackage{url}
\usepackage{xcolor}
\definecolor{newcolor}{rgb}{.8,.349,.1}

\usepackage[linesnumbered]{algorithm2e}
\usepackage{pgfplots,pgfplotstable}
\pgfplotsset{compat=1.14}
\usepackage{tabularx}

\newcommand{\bx}{\ensuremath{\mathbf x}}

\newcommand{\bu}{\ensuremath{\mathbf u}}
\newcommand{\bv}{\ensuremath{\mathbf v}}
\newcommand{\be}{\ensuremath{\mathbf e}\xspace}

\newcommand{\by}{\ensuremath{\mathbf y}}
\newcommand{\bz}{\ensuremath{\mathbf z}}

\newcommand{\bh}{\ensuremath{\mathbf h}}
\newcommand{\bbe}{\ensuremath{\mathbf b}}

\DeclareMathOperator{\tr}{trace}
\DeclareMathOperator{\Var}{Var}
\DeclareMathOperator{\E}{E}

\journal{Journal of Computational Physics}

\begin{document}

\verso{Eloy Romero \textit{etal}}

\begin{frontmatter}

\title{Multigrid deflation for Lattice QCD}

\author[1]{Eloy \snm{Romero}\corref{cor1}}
\cortext[cor1]{Corresponding author}
\ead{eromeroalcalde@wm.edu}
\author[1]{Andreas \snm{Stathopoulos}}
\ead{andreas@cs.wm.edu}
\author[2,3]{Kostas \snm{Orginos}}
\ead{knorgi@wm.edu}

\address[1]{Department of Computer Science, The College of William \& Mary, USA}
\address[2]{Department of Physics, The College of William \& Mary, USA}
\address[3]{Jefferson Laboratory, USA}

\received{10 Sep 2019}
\finalform{--}
\accepted{--}
\availableonline{--}
\communicated{--}

\begin{abstract}
Computing the trace of the inverse of large matrices is typically addressed through statistical methods. Deflating out the lowest eigenvectors or singular vectors of the matrix reduces the variance of the trace estimator.
This work summarizes our efforts to reduce the computational cost of computing the deflation space while achieving the desired variance reduction for Lattice QCD applications. Previous efforts computed the lower part of the singular spectrum of the Dirac operator by using an eigensolver preconditioned with a multigrid linear system solver. Despite the improvement in performance in those applications, as the problem size grows the runtime and storage demands of this approach will eventually dominate the stochastic estimation part of the computation.

In this work, we propose to compute the deflation space in one of the following two ways. First, by using an inexact eigensolver on the Hermitian, but maximally indefinite, operator $A \gamma_5$. Second, by exploiting the fact that the multigrid prolongator for this operator is rich in components toward the lower part of the singular spectrum. We show experimentally that the inexact eigensolver can approximate the lower part of the spectrum even for ill-conditioned operators.
Also, the deflation based on the multigrid prolongator is more efficient to compute and apply, and, despite its limited ability to approximate the fine level spectrum, it obtains similar variance reduction on the trace estimator as deflating with approximate eigenvectors from the fine level operator. 
\end{abstract}


\end{frontmatter}

\section{Introduction}
\label{intro}

To estimate \emph{ab initio} elementary particle properties, Lattice Quantum Chromodynamics (LQCD) involves the evaluation of correlation functions expressed as traces over products of propagators with constant matrices. The propagators include the inverse of the Dirac operator defined on a 4- or 5-dimensional lattice. In practice, the resulting matrix size makes it impossible to compute these traces directly by evaluating the propagators at every point in the lattice. The general strategy is to replace the exact computations by stochastic estimations and control its variance.  
For the case study in this paper, the propagators are closed, which leads to the evaluation of the trace on the inverse of the Dirac operator, $A^{-1}$.

There are several methods to estimate the trace of an inverse matrix. The
simplest method is to estimate the trace as a multiple of the average diagonal element of $A^{-1}$, $\E[N\be_i^\dagger A^{-1}\be_i]$, where $\be_i$ is the $i$-th column of the identity matrix, and $N$ is the dimension of $A$. The variance of this estimator grows with $N^2$ and can be very large even when the diagonal elements of $A^{-1}$ have a small variance.
Sophisticated versions of this approach can be found in \cite{tang2012probing,Wu16}.

The standard tool for this problem is the Hutchinson method \cite{hutchinson} which 
estimates the trace of $A^{-1}$ as 
\begin{equation}
  t(A^{-1}) = \E[\bx^\dagger A\bx] = \frac{1}{s} \sum_{j=1}^{s} \bx_j^\dagger A^{-1}\bx_j,
\end{equation}
where $\bx$ is a random vector whose elements $x_i$ satisfy $\E[x_i^\dagger
x_j] = \delta_{i,j}$.

If the $\bx_i$ are Rademacher vectors, i.e., following the discrete uniform distribution with values
$\pm 1$, the trace estimator can also be written as $t(A^{-1}) = \frac 1 s A^{-1} \odot XX^\dagger$, where $X$ is a matrix with all the $\bx_i$ as columns and $\odot$ is the element-wise product of matrices. Then, the variance associated with the trace estimator $t(A^{-1})$ is
\begin{equation}
  \Var[t(A^{-1})] = \frac{2}{s} \|A^{-1}\|_F^2 - \frac{2}{s}\sum_{i=1}^N |(A^{-1})_{i,i}|^2.
  \label{eq:var}
\end{equation}

In practice $A^{-1}\bx_j$ is computed approximately as $ A^{-1}_\xi(\bx_j) $ by an iterative solver that stops when $\|A\, A^{-1}_\xi(\bx_j) - \bx_j\|\leq \xi\,\|\bx_j\|$. The bias of the resulting operator can be corrected by an independent estimation of the trace of $A^{-1} - A^{-1}_\xi$. This trick is known as the Truncated Solver Method \cite{BALI20101570}.

As in any Monte-Carlo process, the variance reduces linearly with the number of samples, $s$. Every sample involves solving a linear system with the Dirac operator, and for cases where the variance is too large, collecting enough samples can be expensive. Partitioning and deflation can reduce the variance further.

A deterministic way to reduce the variance is by partitioning the matrix rows and columns into $k$ groups and approximating the trace of the diagonal blocks \cite{BERNARDSON1994256,FOLEY2005145,PhysRevD.85.054510}. As $k$ approaches $N$, the variance reduces to zero. In practice, it is used to remove known strongly coupled components on $A^{-1}$ by setting them on different groups. Examples of this technique are spin-color dilution (separating strongly coupled spin-color components) and partitioning (separating close spacetime components). The most basic way of lattice spacetime partitioning is the red-black (immediate neighbors have the opposite color to the node). Hierarchical probing \cite{HP} is a sophisticated extension of this approach that identifies a series of nested colorings at increasing lattice distance. The nesting ensures that, in case of increasing the number of partitions, only linear systems associated with the new colors have to be solved. The drawback of this technique is that the variance is significantly reduced only at the particular number of partitions $k=2^{Di}$, for $i=0,1,\dots$ and $D$ being the number of dimensions of the lattice. Those particular values of $k$ are called the closing colors. 

For similar configurations to the $32^3\times 64$ lattice used in this paper, and with spin-color dilution in all cases, the red-black partition reduces the variance of the original trace estimator by 3-fold, while using higher-level hierarchical probing with 512 vectors reduces the variance by 10-fold \cite{PhysRevD.92.031501,gambhir2017deflation}. The latter reduction in variance is modest compared with merely using Hutchinson with the same number of random vectors. In practice, partitioning is combined with other techniques and especially deflation. 

The initial motivation for deflation is that singular values near zero of the Dirac operator are responsible for the large variance in the trace estimator. Equation~\eqref{eq:var} relates the variance with the spectrum of the traceless matrix inverse through the Frobenius norm, which is the addition of all singular values squared. The singular value spectrum of the traceless matrix is a small perturbation of the original spectrum if the largest singular values of $A^{-1}$ are far greater than the diagonal elements of $A^{-1}$ in magnitude, which is the case for Dirac operators. As shown in \cite{gambhir2017deflation}, 
if the left and the right singular vectors of a non-Hermitian matrix
$A$ are statistically independent standard random unitary matrices, then the
variance of $t(A^{-1})$ is approximately $\sum_{i=1}^N
\sigma_i^{-2}$ (Theorem~2.6 in \cite{gambhir2017deflation}), where $\sigma_i$ are the singular values of $A$.
The authors also show various examples where the variance is close to that estimation,
even when the assumptions of the theorem on the singular vectors are slightly relaxed.

The above suggests that the trace estimator has less variance if applied to an operator removing the near-null part of the spectrum (deflation).
The strategy is, therefore, to split the trace into two parts, $A^{-1}P$ and $A^{-1}(I-P)$, so that the trace of $A^{-1}P$ can be efficiently computed, and the variance of estimating the trace of $A^{-1}(I-P)$ is far smaller than the variance with $A^{-1}$. If the singular values of $A^{-1}$ decay \emph{quickly}, only a small percentage of the spectrum needs to be deflated to reduce sufficiently the variance on the trace estimator of $A^{-1}(I-P)$. For the trace of $A^{-1}$, the left singular vectors of $A$ should be used in $P$ \cite{gambhir2017deflation}.

\begin{figure}
\mbox{
\includegraphics[width=.33\textwidth]{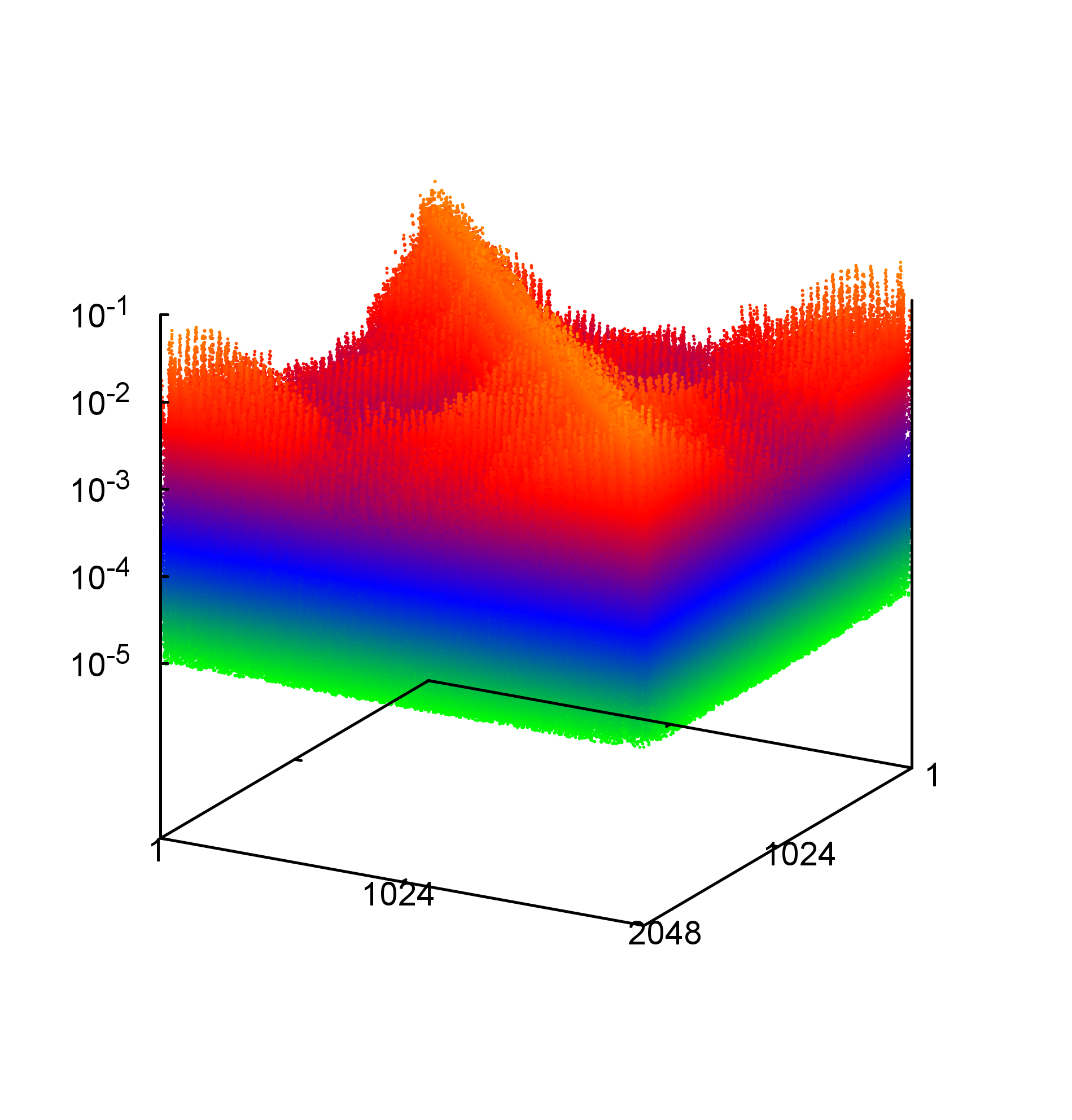}%
\includegraphics[width=.33\textwidth]{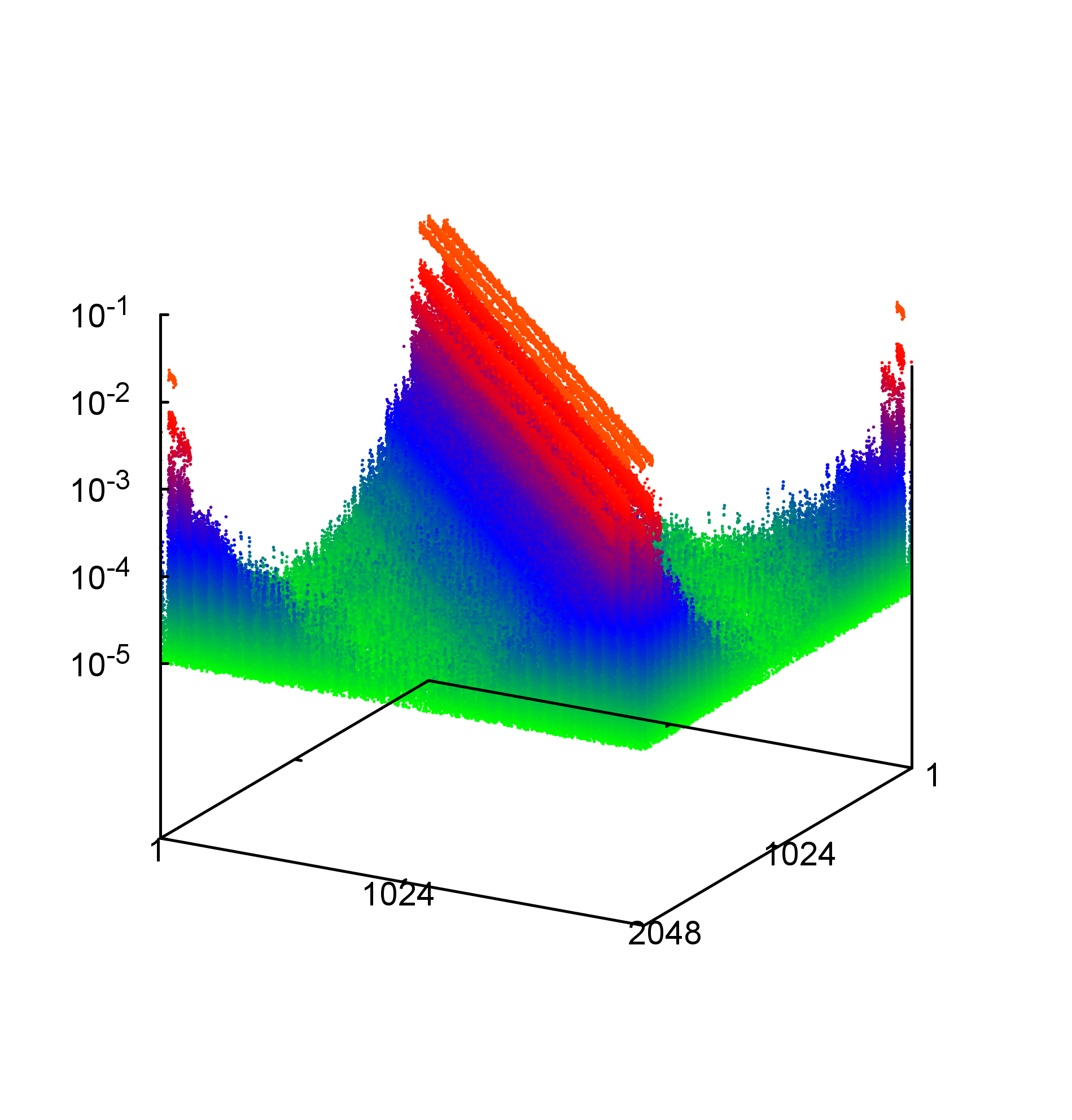}%
\includegraphics[width=.33\textwidth]{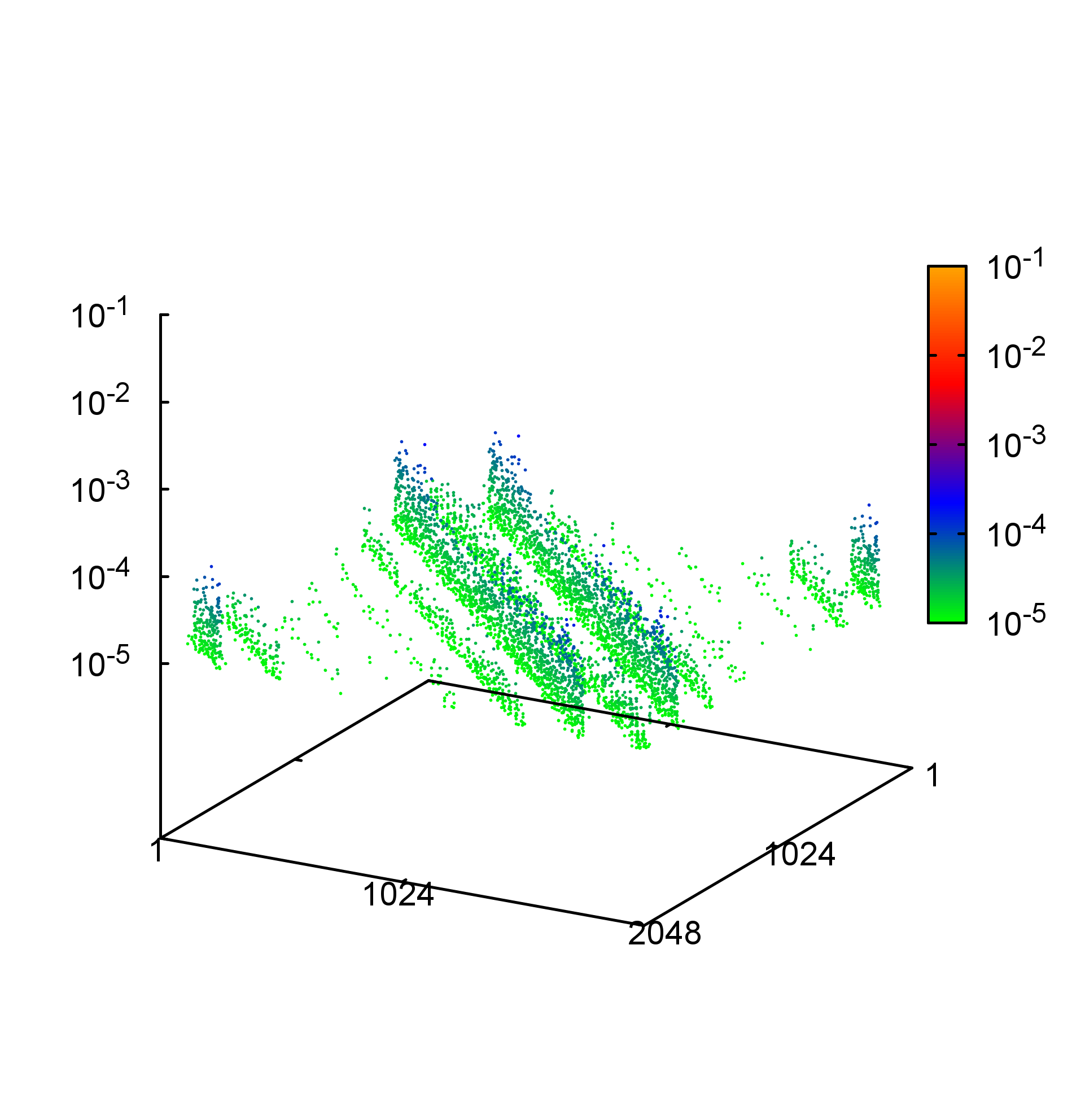}}
\begin{scriptsize}
\newcolumntype{C}{>{\centering}X}
\begin{tabularx}{\textwidth}{CCC}
	 $A^{-1}$ & $A^{-1}(I-UU^\dagger)$ & $(A^{-1}(I-UU^\dagger))\odot (HH^\dagger)$
\end{tabularx}
\end{scriptsize}
\caption{Absolute value of the elements of the inverse of a Wilson 2D $32^2$ (left), after deflating the space associated to the 100 smallest singular values (middle), and, besides that, applying spin-color dilution and 32 basis of hierarchical probing (right).}
\label{fig:deflation}
\end{figure}

An interesting synergy occurs when combining partitioning and deflation. The lower part of the operator's spectrum has the dominant contribution to the magnitude of the elements of $A^{-1}$, mainly affecting the matrix elements that correspond to long distances in the lattice (see Fig.~\ref{fig:deflation}, left). After removing part of the lower spectrum, the remaining eigenvectors contribute to more localized regions of the matrix (see middle graph), for instance, among certain subdiagonals of the matrix, that can be better addressed by partitioning (see right graph). For the same lattice $32^3\times 64$ that hierarchical probing reduces the variance by 10-fold over the original Hutchinson method; also deflating 1000 approximate singular values gives a 100-fold total variance reduction \cite{gambhir2017deflation}. In \cite{BALI20101570}, they reported a speedup of two in cost by combining both techniques. They also considered the acceleration stemming from solving linear systems with the deflated operator, $A^{-1}(I-P)$, which is a well-known effect \cite{eigCG,Abdel2010,AMGLattice}.

To compute the deflation space, previous efforts computed the eigenpairs of the Hermitian Dirac operator, $A \gamma_5$, usually with Chebyshev-accelerated Arnoldi \cite{neff2001low,BALI20101570}.
\cite{gambhir2017deflation} briefly discussed alternatives to compute the deflation space. In particular, they discourage the use of methods such as eigCG \cite{abdel2014extending} or Lanczos because they have to be used unpreconditioned, and thus the convergence is slow. Instead they show results using the state-of-the-art preconditioned eigensolver PRIMME \cite{PRIMME}\footnote{For an introduction to preconditioned eigensolvers see \cite{knyazev1998preconditioned} and section 11 in \cite{etemplates}.} on $A^\dagger A$. The preconditioner employed is based on the adaptive algebraic multigrid (AMG) \cite{AMGLattice} as implemented in QOPQDP\footnote{http://usqcd-software.github.io/qopqdp/}. In this paper, we consider a different strategy for preconditioning that seems more effective.

Because the number of singular vectors required scales as a small percentage of the lattice volume, for large-size lattices the cost of computing and storing the singular vectors becomes unaffordable.
Iterative eigensolvers, e.g., Krylov and Davidson-type, access the operator at least once for every singular vector computed, yielding a quadratic cost in terms of lattice volume\footnote{The block version of the eigensolvers can reduce the number of accesses and the computational time by a constant factor, but it does not remove the quadratic term of the cost.}. Similarly, storage and the application of the deflation basis involve quadratic terms. In addition, orthogonalization of eigenvectors introduces a cubic cost on the number of singular vectors, and thus on the lattice volume, albeit with a much smaller constant because of more efficient memory access patterns.  

In LQCD, the Dirac operator becomes more ill-conditioned as the fermion mass is lighter and, as a consequence, linear system solvers take more iterations (critical slowing down). The adoption of multigrid techniques has accelerated the solution of the linear systems dramatically, not only keeping the number of iterations independent of the operator's conditioning but also making the cost of the solvers scale linearly to the lattice volume \cite{AMGLattice,brannick2008adaptive}.

In this work we propose the use of multigrid directly into reducing the variance of the trace estimator to overcome the scalability limitations of deflation. The approach is described in section~\ref{sec:multigrid}. Although we compute the lower part of the singular spectrum of a much smaller dimension operator, its performance continues to be critical, as we show in section~\ref{sec:computing}. In section~\ref{sec:results}, we evaluate the variance reduction of the new approach experimentally. We offer some conclusions in Section~\ref{sec:conclusions}.

\section{Multigrid Variance reduction}
\label{sec:multigrid}

To establish notation, we show the two-grid method for solving $A\bx=\bbe$ for a Hermitian matrix $A$. A two-grid iterative method alternates between the following two steps:

\begin{enumerate}
\item Coarse grid correction: $\bx_i^c = \bx_i + V(V^\dagger AV)^{-1}V^\dagger (\bbe - A \bx_i)$
\item Smoothing: $\bx_{i+1} = \bx_i^c + M^{-1}(\bbe - A\bx_i^c)$.
\end{enumerate}

In step 1, the method projects the residual vector to a lower dimension (coarsening), corrects it, and prolongates it back. The matrices $V$ and $V^\dagger$ are the multigrid's prolongator (or interpolator) and restrictor respectively.
In step 2, the new residual vector $\bbe-A\bx_i^c$ is \emph{smoothed} by applying few iterations of a linear system solver, e.g., GCR, on $A$ which we represent by the operator $M^{-1}$. 
The goal of the coarse grid is to correct directions that the fine grid smoothing struggles to fix. One of the most relevant distinctions between multigrid variants is the coarsening.

In adaptive multigrid for LQCD, the setup starts with the approximate solutions of $A\by^{(i)}=0$ for $n$ different random initial guesses. The returned $n$ solution vectors, $Y$, are rich in directions of singular vectors with singular values close to the origin, because linear solvers have difficulty resolving these lower spaces. Interestingly, these singular vectors on the lower part of the spectrum look similar on a local lattice scale, and nearby singular vectors can be generated by combining local pieces from different singular vectors \cite{Luscher2007,Clark-Jung-MGdeflation}. Taking advantage of that, the adaptive multigrid method partitions the lattice into a set of $m$ compact subdomains that maintain these ``local'' features. Then, each of the solution vectors is broken into $j=1,\ldots,m$ vectors, each with support only on the subdomain $\mathbb D_j$. This process is referred to as \emph{blocking},
\[
   \left[ Y \right]_D = \left[ D_1 Y,\ \ D_1 Y,\ \ \dots,\ \ D_m Y \right], \text{\ with\ } \mathbb{D}_i \cap \mathbb{D}_j=\emptyset\ \forall i\neq j,
\]
where $D_i$ is the projector onto the domain $\mathbb D_i$. For numerical stability, implementations work with an orthonormal basis $V$ of the blocked solutions $\left[ Y \right]_D$.

To see the role of deflation in multigrid, we rewrite the previous two steps in terms of the oblique projector $P=AV(V^\dagger AV)^{-1}V^\dagger$ on the space of $V$, 
\[
  \bx_{i+1} = \bx_i + A^{-1}P(\bbe - A \bx_i) + M^{-1}(I-P)(\bbe - A \bx_i).
\]
The coarse grid correction is interpreted as the exact solution of the components of the residual vector on $P$, and smoothing takes care of the remaining components of the residual vector.
The deflation view suggests that, in general, if $V$ approximates the lower part of $A$'s spectrum, then the effective condition number of $A^{-1}(I-P)$ is smaller than that of $A^{-1}$, accelerating the convergence of the iterative method $M^{-1}$.

In the case of the trace estimation, the deflation splits the trace of $A^{-1}$ into the trace of $A^{-1}(I-P)$ and the trace of $A^{-1}P$, so that the trace of $A^{-1}P$ can be computed efficiently, and the variance of estimating the trace of $A^{-1}(I-P)$ is far smaller than with $A^{-1}$.
As with solving linear systems, we are interested in $P$ approximating the lower part of $A$'s spectrum.

\subsection{Choice of projector}

If $P$ is the oblique projector defined above and $V$ is of much smaller rank than $N$, the trace of $A^{-1}P$ can be computed efficiently since
$\tr A^{-1} AV(V^\dagger AV)^{-1}V^\dagger = \tr V(V^\dagger AV)^{-1}V^\dagger = 
	\tr (V^\dagger AV)^{-1}$, which is the trace of a smaller matrix.
If $V$ is an invariant subspace, $P$ becomes the orthogonal spectral projector.
If $V$ is far from an invariant space, the variance of the trace estimator $t(A^{-1}(I-P))$ is usually smaller with the orthogonal projector than with the oblique projector. Unfortunately, there is no efficient way to approximate the trace of $A^{-1}P$ if $P=VV^\dagger$, unless $V$ is close enough to left singular space of $A$ so that the errors on the singular values are smaller than the standard deviation of the trace estimator. In this work, we propose the use of the oblique projector but also spend the time to compute and report the variance of $t(A^{-1}(I-P))$ for the orthogonal projector as a rough estimate of the lower bound to the variance achievable by any projector of that rank.

If $A$ is not a Hermitian positive definite matrix, and $V$ is not close to a left singular space, the coarse operator $V^\dagger AV$ may have singular values that are close to zero but do not correspond to singular values of $A$, which are known as spurious values. In addition to making $V^\dagger AV$ ill-conditioned, these small spurious values increase the variance of the $t(A^{-1}(I-P))$. 
One can consider variants of the oblique projector that may prevent the directions in $V$ that are far from an eigenspace from harming the final variance,
\[
  \tr A^{-1} = \tr A^{-1}(I - AV(V^\dagger KAV)^{-1}V^\dagger K) + \tr V^\dagger K V(V^\dagger KAV)^{-1}.
\]

Setting $K=A^{-1}$ avoids the problems of spurious values and makes $P$ the orthogonal projector on $V$.
If $A$ is not Hermitian, another obvious choice for $K$ is $A^\dagger$, but it makes the coarse operator more expensive by introducing more nonzero elements. In the results presented in Section~\ref{sec:results}, we use $K=\gamma_5$, that does not protect against the directions in $V$ that are far from eigenspaces and may result in small spurious eigenvalues, but it is somehow optimal as we explain in the following. For the spurious values, we supplement this technique with filtering of the prolongator space, as explained in the next section.

If $A^{-1}$ is not Hermitian, one can reduce the variance by evaluating the trace of the symmetric operator $(A^{-1} + A^{-T})/2$ instead, which has the same trace and less or equal variance. So symmetric and Hermitian matrices are optimal for variance reduction under \emph{symmetrization}. The spin-color dilution probing basis $H$ transforms any $\gamma_5$-Hermitian matrix into a Hermitian, and we wish that the deflated operator be also $\gamma_5$-Hermitian. It is easy to see that the operator with the projector for $K=\gamma_5$, $A^{-1}(I - AV(V^\dagger \gamma_5 AV)^{-1}V^\dagger \gamma_5 )$, is $\gamma_5$-Hermitian.

\subsection{Filtering}

The blocking of the approximate solutions to the $n$ equations $A\by_i=0$ (the $\by_i$ are often called the \emph{null vectors}) produces approximations to a large part of the lowest singular spectrum of $A$. However, the quality of these approximations deteriorates rapidly beyond a few lowest singular vectors. Still, using these approximations as an orthogonal projector would effectively reduce the variance of the trace estimator. With the oblique projector, there is an optimal size for the deflation space beyond which the interior singular vectors are not accurate at all, and variance increases. We have confirmed this experimentally with Laplacian matrices and the ensembles presented in Section~\ref{sec:results}. 

It is not clear how to identify this optimal size of the deflation basis \emph{a priori}. However, this is not necessary since the application of the oblique deflation operator can be decoupled from the inversions of $A^{-1}\bz_i$ in the Hutchinson method.
Specifically, we propose to compute and store first a sufficiently large rank of the lower part of the spectrum of the coarse operator $V^\dagger A V$, say $\bar U\Sigma \bar V^\dagger$. The deflated operator in the Hutchinson method would then involve 
   $A^{-1}(I-P) \bz_i = A^{-1}\bz_i - V\bar V(\bar U^\dagger V^\dagger K AV\bar V)^{-1}\bar U^\dagger V^\dagger K \bz_i$, 
for the entire $\bar U, \bar V$ or some column subset of them.  
After completing the inversions with the noise vectors $Z = \left\{ \bz_i, i=1,\ldots,n \right\}$, the deflation part can be recomputed with different number of triplets, and the optimal rank of $P$ can be identified.
This \emph{a posteriori} analysis can be done efficiently without storing any vector of lattice dimension, just the matrices, $Z^\dagger V\bar V$, and $\bar U^\dagger V^\dagger KZ$.

Finally, we note the difference between our approach, which computes the eigenspaces from the coarse space,  and the approach in \cite{Clark-Jung-MGdeflation} where the computed fine grid eigenspace is stored compactly using the prolongator. 

\section{Computation of the deflation space}
\label{sec:computing}

The $\gamma_5$-Hermiticity of the Wilson operator as well as of the coarse operators built by multigrid allows for several algorithmic alternatives for computing the smallest singular values and their vectors. First, we examine the relative performance of these alternatives. However, as the density of singular values near zero increases with lattice volume, the eigenvalue problems become particularly challenging for iterative methods, and complementary techniques must be employed.

\subsection{Choice of eigensolver}

A common approach to approximate the smallest singular values of $A$ is to use some variant of the Lanczos method on either $A^\dagger A$ or $A \gamma_5$. 
However, without preconditioning, Lanczos converges slowly for the operators that we are working with.
Alternatively, the shift-and-invert Lanczos has much faster convergence at the expense of solving a linear system of equations per (outer) iteration. 
We consider the inexact Lanczos variant, in which the linear systems are solved approximately with a relative tolerance $\xi$.
The approximation also limits the accuracy of the computed eigenspace by $\xi$.

Besides inexact Lanczos, we also consider Generalized Davidson, a preconditioned method that offers a few more techniques for improving convergence and robustness (the GD+k method) \cite{JDQMR_One}. We test both the inexact GD+k variant on the operator $A^{-1}_{\xi}$ and the exact GD+k on the operator $A\gamma_5$ or $A^\dagger A$, while providing the solution of linear systems as a preconditioner instead.

\begin{table}
	\caption{Number of inversions in computing the 10 smallest singular vectors on a Schwinger model Wilson $32^2$ with several approaches. The eigensolvers mark an eigenpair $(\lambda,\bx)$ as converged when $\|\gamma_5 A^{-1}\bx - \lambda^{-1}\bx\| \leq 10^{-2}\|A^{-1}\|$. The inverter used in all cases, $A^{-1}_\xi$, is QMR stopping at relative tolerance $\xi=10^{-2}$. GD+$k$ is from PRIMME, and Lanczos is Matlab's \texttt{eigs}, both with default settings.}
	\label{fig:deflation2d}

	\begin{center}
	\begin{footnotesize}
	\begin{tabular}{lr}
	Approach & \# Inversions\\\hline
	GD+$k$ on $A^\dagger A$ with preconditioner $A^{-1}_\xi {A^{-1}_\xi}^\dagger$ & 60\\
	GD+$k$ on $A^{-1}_\xi {A^{-1}_\xi}^\dagger$ & 34\\
	GD+$k$ on $A \gamma_5$ with preconditioner $\gamma_5 A^{-1}_\xi $ & 27\\
	GD+$k$ on $\gamma_5 A^{-1}_\xi$ & 23\\
	Lanczos on $A^{-1}_\xi {A^{-1}_\xi}^\dagger$ & 100 \\
	Lanczos on $\gamma_5 A^{-1}_\xi $ & 25
	\end{tabular}
	\end{footnotesize}
	\end{center}
\end{table}

Table~\ref{fig:deflation2d} shows a comparison of the various eigensolver alternatives on a Schwinger model Wilson operator\footnote{Generated with \texttt{quantum-mg}, source code at \url{https://github.com/weinbe2/quantum-mg}.}, a $32^2$ lattice, and $\beta = 6$. The experiments are performed in MATLAB and, although the problem size is small, they are indicative of the performance of these methods on much larger problem sizes. 
We observe that working with $A^\dagger A$ involves more solutions of linear systems than working with $A \gamma_5$. For the latter operator, inexact Generalized Davidson performs slightly better than Generalized Davidson with preconditioner and inexact Lanczos.

\subsection{Eigensolver tuning}

The Dirac operators that we work with exhibit a fast decay of the singular value spectrum, so most of the eigenmethods on $\gamma_5 A^{-1}_\xi$ obtain good approximations within an average of two iterations per singular triplet.
The singular triplets closest to the origin need more inversions (outer steps), and each of those inversions takes many more iterations than for more interior triplets. Therefore tuning the eigensolver and the multigrid inverter can be performed early before computing thousands of singular triplets.

Because $\gamma_5 A^{-1}_\xi$ and the corresponding coarse grid operator are indefinite Hermitian matrices, they still may suffer from many spurious eigenvalues near zero. In some cases, this may cause the inverter to stagnate at an error level well above the requested $\xi$. A common way to avoid stagnation is to shift the operator of the linear systems to reduce its condition number. Shifting the operator by the diagonal $i\gamma_5$ moves the eigenvalues away from the origin (see Fig~\ref{fig:wilson2d}, left), while the eigenvectors of $\gamma_5 A^{-1}$ and $(A\gamma_5 + i\tau I)^{-1}$ have the same order with respect to the distance of their eigenvalues from the origin (see Fig~\ref{fig:wilson2d}, right).

The new operator after the shifting is no longer Hermitian but is normal. One option is to modify a Hermitian eigensolver to support normal operators. Specifically, the projected Rayleigh-Ritz problem at every iteration needs to be solved with a Schur decomposition instead of a Hermitian eigensolver. Another option is to use non-Hermitian eigensolvers.
We followed the former approach, making the proper changes in PRIMME.
When the inverter fails, it happens in the first few iterations of the eigensolver when spurious eigenvalues near zero are present. In that case, we rebuild the multigrid inverter with a larger shift, and the eigensolver is restarted.

A challenge with any eigensolver that computes a large number of eigenpairs is the cost of orthogonalization. The number of floating-point operations associated with orthogonalization grows as 
$O(N\times\text{\#SV}\times\text{iterations})$. In addition, on massively parallel environments, orthogonalization involves between $O(\text{\#SV})$ and $O(\text{\#SV}^2)$ global synchronizations depending on the algorithm. It is clear, therefore, that for sufficiently many singular values ($\text{\#SV}$) orthogonalization time eventually dominates the eigensolver performance for LQCD. 
We demonstrate the performance of the GD+k eigensolver on the configurations of two large lattices shown in Table~\ref{tab:confs}.

\begin{table}
\caption{Configuration of the Dirac-Wilson operator and the multigrid solver used in the experiments. The Preconditioned settings build the multigrid for the odd-even operator for solving linear systems. The Unpreconditioned multigrid is used to obtain the singular vectors of the original operator.}
\label{tab:confs}
\centering
\begin{footnotesize}
\begin{tabular}{cr|cc|llr}
              &   & \multicolumn{2}{c|}{Preconditioned} &  \multicolumn{3}{c}{Unpreconditioned} \\
Configuration & Dimension & Blocking & Null vectors & Blocking & Null vectors & Dimension \\
              & fine grid (N) &  per level (m) & per level (n) & per level & per level & of coarsest
 \\\hline\\[-1.5mm]
$32^3\times 64$                               & 25,165,824  & $4^4,\ 2^4$     & 24, 32  & $4^4$       & 24 & 393,216\\
$\beta=6.3$                                   &             &                 &        & $4^4,\ 2^4$ & 24, 32 & 32,768\\
$\mu=-0.2390$                                 &             &                 &        & $4^4$       & 48 & 786,432\\
                                              &             &                 &        & $4^4,\ 2^4$ & 48, 64 & 65,536\\\hline\\[-1.5mm]
$64^3\times 128$                              & 402,653,184 & $4^4,\ 2^3\times 4$& 24, 32  & $4^4$       & 48 & 12,582,912\\
$\beta=6.3$                                   &             &                 &        & $4^4,\ 2^4$ & 48, 64 & 1,048,576\\
$\mu=-0.2416$ &&&&&&\\
\end{tabular}
\end{footnotesize}
\end{table}

\begin{table}
\caption{Running times for computing the partial singular decomposition on Wilson operators and coarse operators listed on Table~\ref{tab:confs}.}
\label{tab:eigencosts}
\centering
\footnotesize
\begin{tabular}{rrrrr}
Configuration & \# Inversions & Inverter time & Ortho. time & Total time \\\hline\\[-1.5mm]
  \multicolumn{1}{@{}l}{\textit{Wilson $32^3\times 64$ on 22 nodes:}}\\[1mm]
	1024 SV on $A$ $\tau=0.00$                    &  1,429 &   655 &   254 &   957 \\
	nv 24 \& 1024 SV on coarse op.    $\tau=0.04$ &  7,544 & 1,770 &    16 & 1,797 \\
	nv 24, 32 \& 1024 SV on coarse op. $\tau=0.04$ &  6,664 &    90 &     2 &    97 \\
	nv 48 \& 1024 SV on coarse op.    $\tau=0.04$ &  7,794 & 8,570 &    92 & 8,826 \\
	nv 48, 64 \& 1024 SV on coarse op. $\tau=0.04$ &  7,725 &   193 &    22 &   364 \\[1mm]
  \multicolumn{1}{@{}l}{\textit{Wilson $64^3\times 128$ on 256 nodes:}}\\[1mm]
	1024 SV on $A$ $\tau=0.02$                    &  3,821 & 5,874 &   860 & 6,952 \\
	nv 48, 64 \& 1024 SV on coarse op. $\tau=0.02$ &  6,675 &   586 &   119 &   807 \\
	nv 48, 64 \& 4096 SV on coarse op. $\tau=0.02$ & 21,865 & 1,919 & 1,036 & 3,122
\end{tabular}
\end{table}

In Table~\ref{tab:eigencosts} we can see that orthogonalization time starts to dominate for the coarse operator of the Wilson $64^3\times 128$ case. Going from 1024 to 4096 singular triplets, the number and the cost of inversions increases by a factor of 3.3. However, the orthogonalization time increases by a factor of 8.8, making it a third of the total computation time. 
Block variants of iterative methods, including GD+k, can utilize better optimized linear algebra library functions and thus may reduce the orthogonalization time by a small factor. However, orthogonalization remains the scalability bottleneck.

The motivation to find the deflation space from a coarse operator is twofold: compute a much larger space and do all the operations on matrices of much smaller dimension. However, the benefits are not commensurate. First, the coarse level operators $V^\dagger AV$ are typically much denser in nonzeros per row than $A$, so despite a significant reduction in the dimension of the matrix, the time savings are moderate. 
For example, consider the $32^3\times 64$ lattice in Table~\ref{tab:confs}. The first coarse level (24 null vectors) is 12 times smaller than the fine level. Yet, the cost per inversion is only half. 
Second, although the floating-point operations in orthogonalization reduce linearly with the dimension, the number of processors and synchronizations remains the same. Therefore the parallel performance of the orthogonalization kernel deteriorates and its time  does not reduce. In the examples of Table~\ref{tab:confs}, the matrix of the coarsest levels (with 48, 64 null vectors) is 384 times smaller than the fine level. However, when solving for 1024 singular vectors, the orthogonalization times per iteration reduce only by factors of 62 and 12 for the lattices $32^3\times 64$ and $64\times 128$ respectively.

To conclude, even though the inverter on $A$ spends most of the time on the coarsest level, the eigensolver on the coarse operator spends significantly less time on the inverter.

\pgfplotstableread{spectrum_wilson2d.txt}\specwilson
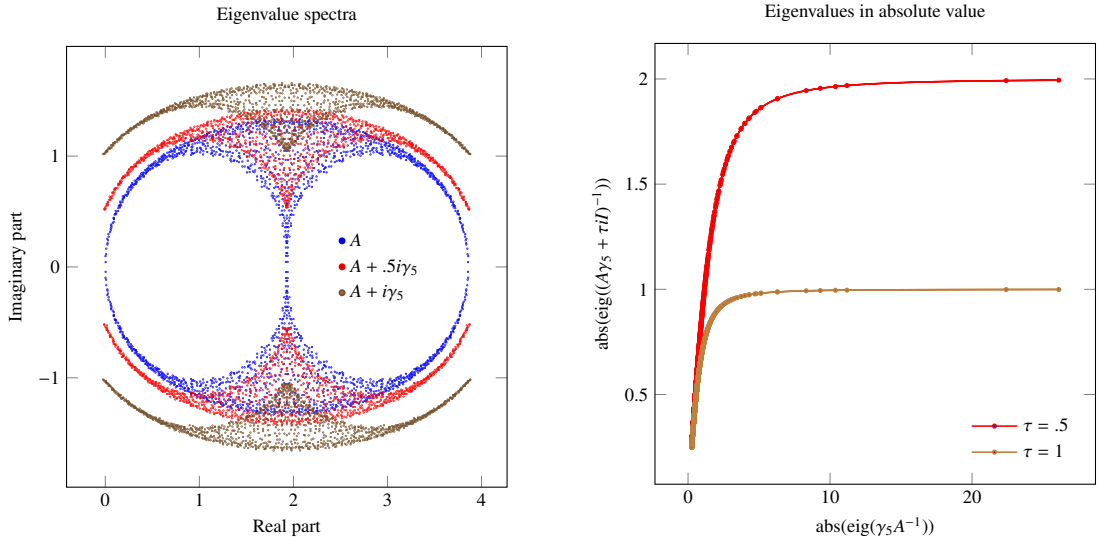
\begin{figure}[t]
\begin{scriptsize}
\hfill
\begin{tikzpicture}[anchor=base]
\begin{axis}[width=.45\textwidth, align={center}, tick label style={/pgf/number format/fixed}, 
        height=.45\textwidth,
        only marks, mark size=0.1pt,
        title = Eigenvalue spectra,
        xlabel = Real part,
        ylabel = Imaginary part,
        every axis plot/.append style={semithick},
        legend entries={$A$,$A + .5i \gamma_5$,$A + i\gamma_5$},
        legend style={draw=none,fill=none,cells={anchor=west, mark size=1pt},anchor=west,at={(.6,.5)}}
  ]
  \addplot plot[mark=*] table[x=eigA-real, y=eigA-imag] {\specwilson};
  \addplot plot[mark=*] table[x=eigAplusip5g-real, y=eigAplusip5g-imag] {\specwilson};
  \addplot plot[mark=*] table[x=eigAplusig-real, y=eigAplusig-imag] {\specwilson};
\end{axis}
\end{tikzpicture}
\hfill
\begin{tikzpicture}[anchor=base]
\begin{axis}[width=.45\textwidth, align={center}, tick label style={/pgf/number format/fixed}, 
        height=.45\textwidth,
        title = Eigenvalues in absolute value,
        xlabel = abs(eig($\gamma_5 A^{-1}$)),
        ylabel = abs(eig($(A\gamma_5 + \tau i I)^{-1}$)),
        every axis plot/.append style={semithick},
        mark size=.65pt,
        legend image post style={mark=*},
        legend entries={$\tau=.5$,$\tau=1$},
        legend pos=south east,
        legend style={draw=none,fill=none,cells={anchor=west}}
  ]
  \addplot plot[no markers,color=red] table[x=abseigginvAint, y=abseiginvAplusigp5int] {\specwilson};
  \addplot plot[mark=*,only marks,forget plot,color=red] table[x=abseigginvA, y=abseiginvAplusigp5] {\specwilson};
  \addplot plot[no markers,color=brown] table[x=abseigginvAint, y=abseiginvAplusigint] {\specwilson};
  \addplot plot[mark=*,only marks, forget plot,color=brown] table[x=abseigginvA, y=abseiginvAplusig] {\specwilson};
\end{axis}
\end{tikzpicture}
\hfill
\end{scriptsize}
\caption{Eigenvalue spectrum after shifting a Schwinger model Wilson $32^2$, $A$, with $\tau i \gamma_5$ for values of $\tau$ of $0.5$ and $1$ (left). Transformation of the eigenvalue spectrum of the inverse of the Hermitian operator, $\gamma_5 A^{-1}$, after shifting, $(A \gamma_5+ \tau i I)^{-1}$ (right).}
\label{fig:wilson2d}
\end{figure}

 


\section{Experimental Performance of the Estimator}
\label{sec:results}

Our benchmark deflation test is the use of the orthogonal projector of the approximate left singular vectors of $A$, computed with an inexact eigensolver as described in Section~\ref{sec:computing} and Algorithm \ref{alg:inexact} below. 
Our proposed method uses an oblique projector to deflate the subspace of a multigrid prolongator that corresponds to the lower part of the singular spectrum of the coarse operator as described in Algorithm~\ref{alg:multigrid}.

\begin{algorithm}
\DontPrintSemicolon
\SetKwInOut{Input}{Input}
\Input{Operator, $A$; deflation space rank, $k$; tolerance, $\xi$; shift, $\tau$; partition basis, $H$}
Compute deflation space: $U = \text{eigs}(\gamma_5(A + i \tau \gamma_5)^{-1}_\xi, \text{largest magnitude}, k, \xi)$\;
Compute the deterministic part: $t_0 = \sum_i \bu_i^\dagger A^{-1}_\xi\bu_i$ \;
Compute the stochastic part: $t_1 = t(A^{-1}_\xi(I-UU^\dagger) \odot HH^\dagger)$ \;
Return $t_0 + t_1$ \;
\caption{Variance reduction for trace estimation using an orthogonal projector of approximate singular vectors of $A$}
\label{alg:inexact}
\end{algorithm}

\begin{algorithm}
\DontPrintSemicolon
\SetKwInOut{Input}{Input}
\Input{Operator, $A$; deflation space rank, $k$; tolerance, $\xi$; shift, $\tau$; partition basis, $H$}
Compute multigrid prolongator $W$; $\gamma_5 W = W\gamma_5^W$ \;
Compute deflation: $\bar U = \text{eigs}(\gamma_5^W (W^\dagger AW+i\tau \gamma_5^W)^{-1}_\xi, \text{largest magnitude}, k, \xi)$ \;
Factorize $\hat U \Lambda\hat U^\dagger = \bar U^\dagger \gamma_5^W W^\dagger A W\bar U$, and set $U=\bar U\hat U$\;
Compute the direct part: $t_0 = \sum_i \bu_i^\dagger \gamma_5\bu_i \lambda_i^{-1}$ \;
Compute the stochastic part: $t_1 = t(A_\xi^{-1}(I-AW\gamma_5U\Lambda^{-1} U^\dagger W^\dagger) \odot HH^\dagger)$\;
Return $t_0 + t_1$ \;
\caption{Variance reduction for trace estimation using an oblique projector from approximate singular subspace of the coarse operator}
\label{alg:multigrid}
\end{algorithm}

The implementation is developed within QOPQDP, a library in the USQCD suite which provides several operators and linear system solvers. We only use the Wilson clover operator and the multigrid even-odd solver. We instantiate two multigrid variants as explained in the caption of Table~\ref{tab:confs}. The preconditioned multigrid variant is based on the odd-even operator and is used for solving the linear systems in the stochastic approximation of the trace, and also for the inversions of the inexact eigensolver of Algorithm~1. The unpreconditioned multigrid variant is used for generating the prolongators and the coarse operators from which the deflation space will be computed in Algorithm~2. At coarse levels, the solution of the linear systems is not accelerated with a multigrid recursion. Instead, even-odd preconditioned GCR is employed. We extended the implementation of QOPQDP solvers to support the implicit $i\tau\gamma_5$ shifting of the operator and the access to the coarse operators from PRIMME.

For the stochastic estimation step, $t(A_\xi^{-1}(I-P) \odot HH^\dagger)$, our implementation follows a blocked approach in which the 12 vectors of spin-color dilution of each hierarchical probing vector $\bh_i$ are processed together. This can potentially accelerate the solution of the linear systems, and the application of the projectors for Algorithm~\ref{alg:inexact}, where the orthogonal projectors can employ optimized BLAS matrix-matrix multiplication routines.

Tuning the prolongators configuration is crucial to get good variance reduction. In the examples tested, we found that the number of null vectors used for setting up the multigrid levels should be twice as many for Algorithm~\ref{alg:multigrid} than the ones used only to solve linear systems. For variance reduction, we configure the solver used on the null vectors to stop at a relative residual norm tolerance of $10^{-2}$ or after 200 iterations.

\subsection{Variance Reduction}

\pgfplotstableread{cl21_32_64_b6p3_m0p2390_m0p2050_cfg_5000-inv24_32_its5086_ig-trace.txt}\clthirtytwosixtyfourinvtwentyfourthirtytwo
\pgfplotstableread{cl21_32_64_b6p3_m0p2390_m0p2050_cfg_5000-inv48_64_its9859_ig-trace.txt}\clthirtytwosixtyfourinvfortyeightsixtyfour
\pgfplotstableread{cl21_32_64_b6p3_m0p2390_m0p2050_cfg_5000-ip.txt}\clthirtytwosixtyfourip
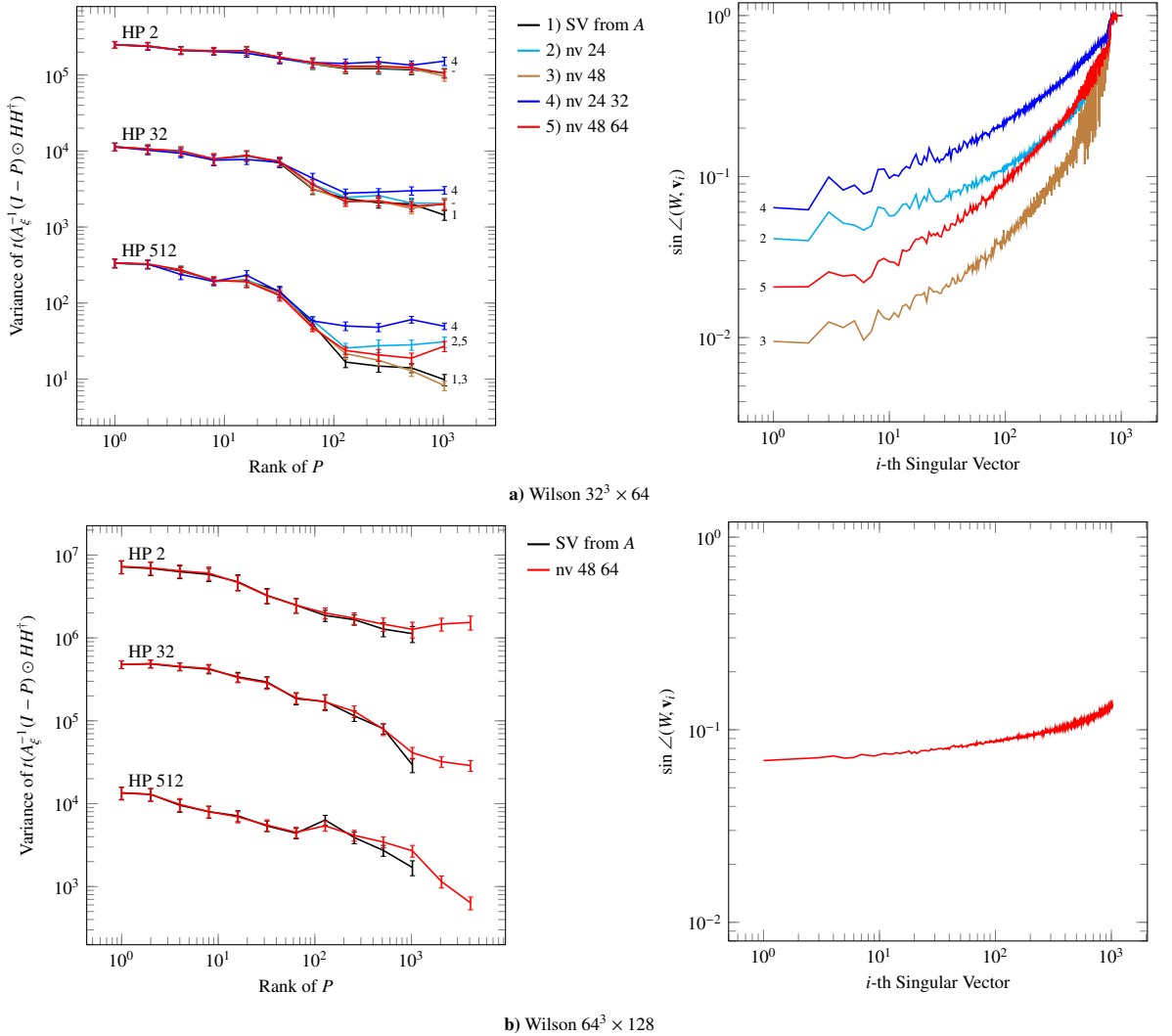
\begin{figure}[t]
\begin{center}
\begin{scriptsize}
\begin{tikzpicture}
\begin{loglogaxis}[width=.44\textwidth, align={center},tick label style={/pgf/number format/fixed}, 
        height=.44\textwidth,
        xmax = 3000,
        no markers,
        xlabel = Rank of $P$,
        ylabel = Variance of $t(A^{-1}_\xi(I-P)\odot HH^\dagger)$,
        cycle list={black,cyan,brown,blue,red},
        every axis plot/.append style={semithick,error bars/y dir=both,error bars/y explicit, error bars/error bar style={solid, semithick, mark size=1pt}},
        legend entries={1) SV from $A$,2) nv 24,3) nv 48,4) nv 24 32,5) nv 48 64},
	legend image post style={xscale=0.5,thick},
        legend style={draw=none,fill=none,cells={anchor=west},legend pos=outer north east}
  ]
  \addplot table[x=eigs, y=level0-hp2-var, y error=level0-hp2-err] {\clthirtytwosixtyfourinvtwentyfourthirtytwo} node[pos=0,anchor=south west,color=black] {HP 2};
  \addplot table[x=eigs, y=level1-hp2-var, y error=level1-hp2-err] {\clthirtytwosixtyfourinvtwentyfourthirtytwo} node[pos=1,anchor=west,color=black] {\tiny -};
  \addplot table[x=eigs, y=level1-hp2-var, y error=level1-hp2-err] {\clthirtytwosixtyfourinvfortyeightsixtyfour};
  \addplot table[x=eigs, y=level2-hp2-var, y error=level2-hp2-err] {\clthirtytwosixtyfourinvtwentyfourthirtytwo} node[pos=1,anchor=west,color=black] {\tiny 4};
  \addplot table[x=eigs, y=level2-hp2-var, y error=level2-hp2-err] {\clthirtytwosixtyfourinvfortyeightsixtyfour};
  \addplot table[x=eigs, y=level0-hp32-var, y error=level0-hp32-err] {\clthirtytwosixtyfourinvtwentyfourthirtytwo} node[pos=0,anchor=south west,color=black] {HP 32} node[pos=1,anchor=west,color=black] {\tiny 1};
  \addplot table[x=eigs, y=level1-hp32-var, y error=level1-hp32-err] {\clthirtytwosixtyfourinvtwentyfourthirtytwo};
  \addplot table[x=eigs, y=level1-hp32-var, y error=level1-hp32-err] {\clthirtytwosixtyfourinvfortyeightsixtyfour};
  \addplot table[x=eigs, y=level2-hp32-var, y error=level2-hp32-err] {\clthirtytwosixtyfourinvtwentyfourthirtytwo} node[pos=1,anchor=west,color=black] {\tiny 4};
  \addplot table[x=eigs, y=level2-hp32-var, y error=level2-hp32-err] {\clthirtytwosixtyfourinvfortyeightsixtyfour} node[pos=1,anchor=west,color=black] {\tiny -};
  \addplot table[x=eigs, y=level0-hp512-var, y error=level0-hp512-err] {\clthirtytwosixtyfourinvtwentyfourthirtytwo} node[pos=0,anchor=south west,color=black] {HP 512} node[pos=1,anchor=west,color=black] {\tiny 1,3};
  \addplot table[x=eigs, y=level1-hp512-var, y error=level1-hp512-err] {\clthirtytwosixtyfourinvtwentyfourthirtytwo} node[pos=1,anchor=west,color=black] {\tiny 2,5};
  \addplot table[x=eigs, y=level1-hp512-var, y error=level1-hp512-err] {\clthirtytwosixtyfourinvfortyeightsixtyfour};
  \addplot table[x=eigs, y=level2-hp512-var, y error=level2-hp512-err] {\clthirtytwosixtyfourinvtwentyfourthirtytwo} node[pos=1,anchor=west,color=black] {\tiny 4};
  \addplot table[x=eigs, y=level2-hp512-var, y error=level2-hp512-err] {\clthirtytwosixtyfourinvfortyeightsixtyfour};
\end{loglogaxis}
\end{tikzpicture}
\begin{tikzpicture}
\begin{loglogaxis}[width=.44\textwidth, align={center}, tick label style={/pgf/number format/fixed}, 
        height=.44\textwidth,
        no markers,
        xlabel = $i$-th Singular Vector,
        ylabel = {$\sin \angle (W,\bv_i)$},
        cycle list={cyan,brown,blue,red},
        every axis plot/.append style={semithick},
        ymin = 0.003,
        ymax = 1.2,
  ]
  \addplot table[x expr=\coordindex, y=inv24_32_its5086_ig-level1-norm] {\clthirtytwosixtyfourip} node[pos=0,anchor=east,color=black] {\tiny 2};
  \addplot table[x expr=\coordindex, y=inv48_64_its9859_ig-level1-norm] {\clthirtytwosixtyfourip} node[pos=0,anchor=east,color=black] {\tiny 3};
  \addplot table[x expr=\coordindex, y=inv24_32_its5086_ig-level2-norm] {\clthirtytwosixtyfourip} node[pos=0,anchor=east,color=black] {\tiny 4};
  \addplot table[x expr=\coordindex, y=inv48_64_its9859_ig-level2-norm] {\clthirtytwosixtyfourip} node[pos=0,anchor=east,color=black] {\tiny 5};
\end{loglogaxis}
\end{tikzpicture}\\
\textbf{a)} Wilson $32^3\times 64$\\
\medskip
\begin{tikzpicture}
\begin{loglogaxis}[width=.44\textwidth, align={center}, tick label style={/pgf/number format/fixed}, 
        height=.44\textwidth,
        no markers,
        xlabel = Rank of $P$,
        ylabel = Variance of $t(A^{-1}_\xi(I-P)\odot HH^\dagger)$,
        cycle list={black,red},
        every axis plot/.append style={semithick},
        every axis plot/.append style={semithick,error bars/y dir=both,error bars/y explicit, error bars/error bar style={solid, semithick, mark size=1pt}},
        legend entries={SV from $A$,nv 48 64},
	legend image post style={xscale=0.5,thick},
        legend pos=south west,
        legend style={draw=none,fill=none,cells={anchor=west},legend pos=outer north east}
  ]
  \pgfplotsinvokeforeach{2,32,512}{
      \addplot table[x=eigs, y=level0-hp#1-var, y error=level0-hp#1-err] {cl21_64_128_b6p3_m0p2416_m0p2050-800_cfg_1342-inv48_64_ig-trace.txt};
      \addplot table[x=eigs, y=level2-hp#1-var, y error=level2-hp#1-err] {cl21_64_128_b6p3_m0p2416_m0p2050-800_cfg_1342-inv48_64_ig-trace.txt} node[pos=0,anchor=south west,color=black] {HP #1};
  }
\end{loglogaxis}
\end{tikzpicture}
\begin{tikzpicture}
\begin{loglogaxis}[width=.44\textwidth, align={center}, tick label style={/pgf/number format/fixed}, 
        height=.44\textwidth,
        no markers,
        xlabel = $i$-th Singular Vector,
        ylabel = {$\sin \angle (W,\bv_i)$},
        cycle list={red},
        every axis plot/.append style={semithick},
        ymin = 0.008,
        ymax = 1.2,
  ]
  \foreach \level in {2}
    \addplot table[x expr=\coordindex, y=inv48_64_ig-level\level-norm] {cl21_64_128_b6p3_m0p2416_m0p2050-800_cfg_1342-ip.txt};
\end{loglogaxis}
\end{tikzpicture}\\
\textbf{b)} Wilson $64^3\times 128$
\end{scriptsize}
\end{center}
\caption{Experimental estimation of the variance of $t(A^{-1}_\xi(I-P) \odot HH^\dagger)$ as the rank of the deflation projector $P$ increases (left), and the sine of the angles between the approximate singular vectors from the operator and the subset of the prolongators used for deflation (right), for Wilson $32^3\times 64$ (top) and Wilson $64^3\times 128$ (bottom). Error bars indicate the jackknife error. The tolerance is $\xi=10^{-2}$.}
\label{fig:var32}
\end{figure}

We show that deflating with a prolongator subspace reduces the variance of the trace estimator, and the effectiveness is related to how close the prolongator subspace is to the lower part of the operator's spectrum. As a base case for comparison, we consider the undeflated matrix after 2 hierarchical probing (HP) vectors have been applied with full spin-color dilution, i.e., 24 inversions. The best results are obtained after 512 HP vectors with the maximum number of vectors deflated. 

In Fig.~\ref{fig:var32}(a) we show results on the Wilson $32^3\times 64$ lattice. Deflating the 1024 smallest singular triplets computed with Algorithm~\ref{alg:inexact} with tolerance $\xi=10^{-2}$ and  512 HP, spin-color diluted vectors
reduces the variance more than four orders of magnitude. 
The top-left graph in Fig.~\ref{fig:var32}(a) shows the effect of prolongators with the four different configurations listed in Table~\ref{tab:confs}. 
For the HP2 and HP32 cases, deflation does not have as much of an effect because the nearest neighbor connections that HP has not removed yet dominate.
For 512 HP vectors, we see that deflation from the first coarse level is very close to the benchmark fine grid deflation (curves 1 and 3). Deflating with singular spaces from the second level is also close to the benchmark up to 100 singular vectors but deflating with more does not have an additional benefit. 

An explanation for the limited benefit at the coarsest level is seen in the top-right graph of Fig.~\ref{fig:var32} that shows the sines of the angles between the computed singular vectors of $A$ and the prolongator subspaces. The prolongators with smaller angles tend to perform better. The graph also shows how fast these angles increase in the interior of the spectrum. All prolongator subspaces have almost no component to singular spaces after the 1000th lowest one, and we should not expect a further reduction on the variance estimation by increasing the rank.
Finally, notice that prolongators generated with more null vectors (curves 3 and 5) yield a smaller variance.

In Fig.~\ref{fig:var32}(b) we show a similar experiment on the Wilson $64^3\times 128$ lattice. We limit the number of computed triplets from the original operator to 1024 due to memory limitations on Edison's nodes, but we let the multigrid deflation use up to a rank of 4096. As in the previous case, both deflation spaces have a similar effect on the variance for the same rank. Deflation achieves a speedup of 20 over the undeflated HP512 experiment.

The second operator is 16 times larger than the first, so it expected that many more vectors are needed to get a similar deflation effect as on the first operator. At the same time, the prolongator subspace from the second coarse level of the second operator has more components on the inner part of the spectrum. This fact may indicate that the effectiveness of the multigrid deflation also scales with the operator dimensions. Further experiments are required to confirm these claims.

Finally, we give a rough estimate on the speedup of these techniques over performing simply a Hutchinson method with random vectors and no HP or dilution.
Both experiments start at 24 inversions and obtain a little more than 4 orders of magnitude improvements with 512*12 = 6144 inversions. This corresponds to speedups of about 50-100 over the Hutchinson method.

\subsection{Performance}

\begin{table}
\caption{Times on NERSC's Edison estimating the trace from 2 noise vectors ($s=2$) solving the linear systems with a tolerance $\xi=10^{-2}$ using a partitioning of full spin-color dilution and 512 hierarchical proving vectors. Approaches detailed in algorithms \ref{alg:inexact} and \ref{alg:multigrid} are compared.}
\label{tab:costs}
\begin{center}
\begin{footnotesize}
\begin{tabular}{rr|ccc|c}
		& & \multicolumn{3}{c|}{Computing $t(A^{-1}_\xi(I- P) \odot HH^\dagger)$} \\
	\multicolumn{2}{c|}{Computing deflating space $P$} & \hspace{2mm} & {\footnotesize Applying P} & {\footnotesize Solving LS} & Variance \\\hline\\[-1.5mm]
  \multicolumn{1}{@{}l}{\textit{Wilson $32^3\times 64$ on 22 nodes:}}\\[1mm]
	1024 SV on $A$ & 957s & & 353s & 2,400s & 10 \\
	nv 24 \& 1024 SV on coarse op. & 269s & & 297s & 2,400s & 35 \\
	nv 24 32 \& 1024 SV on coarse op. & 97s & & 303s & 2,400s & 55 \\
	nv 48 \& 1024 SV on coarse op. & 5489s & & 318s & 2,400s & 8 \\
	nv 48 64 \& 1024 SV on coarse op. & 430s & & 326s & 2,400s & 18 \\[1mm]
  \multicolumn{1}{@{}l}{\textit{Wilson $64^3\times 128$ on 256 nodes:}}\\[1mm]
	1024 SV on $A$ & 6,952s & & 725s & 5,903s & 1,700 \\
	nv 48 64 \& 4096 SV on coarse op.  & 3,122s & & 1,126s & 5,903s & 640
\end{tabular}
\end{footnotesize}
\end{center}
\end{table}

Table~\ref{tab:costs} reports the timings for the computational kernels involved in estimating the trace with just two noise vectors. The most expensive operations are the computation of the deflation space and the inversions in the stochastic estimation of the trace. Notice that applying the deflation for approximate singular vectors from the original matrix is usually cheaper than applying the oblique projector in multigrid deflation, even though the latter prolongator requires less storage. This is because the oblique projector involves a multiplication with the original operator, and because the orthogonal projector is implemented in a block fashion, which is ten times faster than applying the deflation vector by vector.  Using a block version of the oblique projectors, i.e., applying the operator and the prolongators on multiple vectors at once, may also accelerate the computation significantly.

\section{Summary}
\label{sec:conclusions}

In this work, we show that, in the context of reducing variance in the stochastic evaluation of disconnected loops in LQCD, inexact eigensolvers are efficient in approximating the lower part of the singular value spectrum when a fast linear system solver is available. The Dirac operator can be shifted if the problem is too ill-conditioned.
However, computing and storing a fraction of the lattice volume of the singular vectors is computationally prohibitive.

We show that the prolongators generated by adaptive multigrid solvers without even-odd preconditioning have limited accuracy on the lower part of the spectrum, and the accuracy degrades quickly into the inner part of the spectrum.
However, the decay in the accuracy may be slower on larger Dirac operators. In addition, proper tuning of the multigrid configuration can give enough accuracy to make the prolongators suitable for reducing the variance of trace estimators on the ensembles we have tested and with significant computational gains over the fine-grid deflation. As noted before, for large matrices, the speedup over plain Hutchinson method without hierarchical probing can reach two orders of magnitude.

\section*{Acknowledgements}

This research was supported by the Exascale Computing Project (ECP), Project
Number: \mbox{17-SC-20-SC}, a collaborative effort of two DOE organizations---the
Office of Science and the National Nuclear Security
Administration---responsible for the planning and preparation of a capable
exascale ecosystem---including software, applications, hardware, advanced
system engineering, and early testbed platforms---to support the nation's
exascale computing initiative.
K. O. was supported in part by U.S.  DOE grant  \#DE-FG02-04ER41302 and in part by the Jefferson
Science Associates, LLC under U.S. DOE Contract \#DE-AC05-06OR23177.
This work was performed in part using computing
facilities at the College of William and Mary which were provided by
contributions from the National Science Foundation (MRI grant
PHY-1626177), the Commonwealth of Virginia Equipment Trust Fund and
the Office of Naval Research. In addition, this work used resources at
NERSC, a DOE Office of Science User Facility supported by the Office
of Science of the U.S. Department of Energy under Contract
\#DE-AC02-05CH11231, as well as resources of the Oak Ridge Leadership Computing Facility at the Oak Ridge National Laboratory, which is supported by the Office of Science of the U.S. Department of Energy under Contract No. \mbox{\#DE-AC05-00OR22725}.

\bibliographystyle{model1-num-names}
\bibliography{bio}

\end{document}